\magnification=\magstep1
\input amstex
\voffset=-3pc
\documentstyle{amsppt}
\NoBlackBoxes
\def\script U{\Cal U}
\def\script W{\Cal W}
\def\script S{\Cal S}
\def\In{{in}}
\topmatter
\title
Some Automatic Continuity Theorems for Operator Algebras\\
and Centralizers of Pedersen's Ideal
\endtitle
\rightheadtext{Centralizers of Pedersen's ideal}
\leftheadtext{Lawrence G.~Brown}
\author
Lawrence G. Brown\\
Department of Mathematics\\
Purdue University\\
W. Lafayette, IN 47907, USA\\
e-mail: lgb \@ math. purdue.edu
\endauthor
\abstract We prove automatic continuity theorems for ``decomposable"
or``local" linear transformations between certain natural subspaces of
operator algebras.  The transformations involved are not algebra
homomorphisms but often are module homomorphisms.  We show that all left
(respectively quasi-) centralizers of the Pedersen ideal of a $C^{*}$
-algebra $A$ are locally bounded if and only if $A$ has no infinite
dimensional elementary direct summand.  It has previously been shown by
Lazar and Taylor and Phillips that double centralizers of Pedersen's ideal
are always locally bounded.
\endabstract
\endtopmatter
\document
\subheading{\S 1. Introduction and preliminaries.}

In his doctoral dissertation, \cite{27}, N.C. Wong showed that both $Mp$,
where $M$ is a von Neumann algebra and $p$ a projection in $M$, and $A/L$,
where $A$ is a $C^{*}$-algebra and $L$ a closed left ideal, can be
represented as classes of ``admissible"" sections of appropriate fields of
Hilbert spaces, $\{H_{\varphi}: \varphi\in F\}$.  In both cases $F$ is a set
of positive functionals and $H_{\varphi}$ arises from the GNS construction.
Some of Wong's results concern decomposable linear transformations on $Mp,
A/L$, or related spaces of sections, where a linear transformation $T$ is
called decomposable if it arises from a collection $\{T_{\varphi}:
\varphi\in F\}$ of linear transformations on $H_{\varphi}$.  More
precisely, $(Tf)(\varphi)=  T_{\varphi}f(\varphi)$ when $f$ is an
admissible section on $\{H_{\varphi}\}$ (in particular $f(\varphi)\in
H_{\varphi})$.  

Wong's results require that $T$ be bounded, and it is natural to ask
whether decomposable linear transformations are automatically bounded.
This question motivated the present paper, and we show that the answer is
sometimes, but not always, ``yes".  Since our results are in a more general
context than indicated above, we now provide definitions of
decomposability that do not use, or require any knowledge of, fields of
Hilbert spaces or sections.  More explilcit descriptions of the results
will be given after the definitions.

If $x$ is an element of a von Neumann algebra $M, \In \ x$ denotes the
right support projection of $x$; i.e., the support projection of $x^{*}x$.  

\noindent {\it Definition. 1.1.} If $X$ is a subspace of $M$ and $T$ is a linear
transformation from $X$ to $M$, then $T$ is called {\it decomposable} if
$\In \ Tx\leq \In \ x, \forall x\in X$.

The proof of the following trivial result is left to the reader.
\proclaim{Lemma 1.2} $T$ is decomposable if and only if
$\varphi(x^{*}x)=0$ implies $\varphi[(Tx)^{*}(Tx)]=0$ for all normal
states $\varphi$ on $M$ and all $x$ in $X$.
\endproclaim
\noindent {\it Remark 1.3} If $T$ is decomposable and $X \subset Mp$ for
some projection $p$, then $T(X)\subset Mp$.  Also, if $X\subset Mp$ and
$T: X\to Mp$, it is sufficient to verify the condition of 1.2 for $\varphi$
supported by $p$.

Theorem 2.1 below states that if $p$ and $q$ are projections in $M$, then
all decomposable linear transformations from $qMp$ to $Mp$ are bounded if
and only if the following is true:  If $z$ is any central projection of $M$
such that $zpMp$ is finite dimensional, then $zqMp$ is finite dimensional.
Results of Wong imply that any bounded decomposable linear transformation
from $qMp$ to $Mp$ is of the form $x\mapsto ax$, for some $a$ in $Mq$.

If $A$ is a $C^{*}$-algebra, then $A^{**},$ the bidual of $A$, can be
regarded as a von Neumann algebra, usually called the enveloping
$W^{*}$-algebra of $A$.  If $\pi: A\to B(H)$ is the universal
representation of $A$, then $A^{**}$ can be identified with
$\pi(A)^{\prime\prime}$,  the double commutant.  Also, the normal state
space of $A^{**}$ can be identified with the state space of $A$.  The
reader is referred to standard references, for example sections 3.7 and
3.8 of \cite{22}.

Now the concept of decomposability can be applied with $A^{**}$ playing the
role of $M$, but for applications it is better to use a weaker concept.
Let $P(A)$ denote the set of pure state of $A$.

\noindent{{\it Definition 1.4.} If $X$ is a subspace of $A^{**}$ and $T$ is
a linear transformation from $X$ to $A^{**}$, then $T$ is called {\it
purely decomposable} if $\varphi(x^{*}x) = 0$ implies
$\varphi[(Tx)^{*}(Tx)]=0$ for all $\varphi$ in $P(A)$ and all $x$ in $X$.

\noindent{\it Remark1.5.} Again it is true that if $X, T(X)\subset
A^{**}p$ for some projection $p$ in $A^{**}$, then it is enough to check
the condition for pure states supported by $p$.  In the main applications
this will be so for a special kind of projection $p$.  In order to deduce
that $T(X)\subset A^{**}p$ from $X\subset A^{**}p$ and pure
decomposability, one needs extra hypotheses, which will be satisfied in
the main applications.

Certain projections in $A^{**}$ are designated as closed or open, where $p$
is closed if and only if $1-p$ is open; and there is a one-to-one correspondence between closed left ideals $L$ of $A$ and closed projections $p$ in
$A^{**}$, in which $L= \{x\in A: xp=0\}$ (\cite{1},\cite{14}).  Also, $A/L$ can
be isometrically identified with $Ap$, where $x+L$ corresponds to $xp$
(cf. \cite{5, Prop. 4.4}, \cite{8, Thm. 3.3}).  A special case of Theorem
3.6 below is that if $p$ is a closed projection in $A^{**}$ such that there
is no infinite dimensional irreducible representation $\pi$ of $A$ with
$\pi^{**}(p)$ of non-zero finite rank, and if $T$ is a purely decomposable
linear transformation from $Ap$ to itself, then $T$ is bounded.  The full
theorem is stated in a more abstract way.  In the cases of most interest,
if $T$ is bounded and purely decomposable, then $T$ is a left
multiplication, so that in particular $T$ is decomposable (cf. Thm. 5.11 of
\cite{27}).

Every $C^{*}$-algebra $A$ has a smallest dense (two-sided) ideal called the
{\it Pedersen ideal} and denoted by $K(A)$ (see \cite{22, \S 5.6}).  If $A$
is the commutative $C^{*}$-algebra $C_{0}(X)$ for a locally compact
Hausdorff space $X$, then $K(A)$ is the set of continuous functions of
compact support.  Lazar and Taylor \cite{18} and Phillips \cite{25} showed
that any double centralizer $(T_{1}, T_{2})$ of $K(A)$ is locally bounded
in the following sense: There is a familly $\{I_{j}\}$ of closed two-sided
ideals of $A$ such that $(\Sigma I_{j})^{-} = A$ and $T_{1}$ and $T_{2}$
are bounded on $K(A)\cap I_{j}$ for each $j$.  The analogous statements
for left or quasi-centralizers are false in general, as was pointed out to
us by H. Kim, but in Section 4, Theorem 3.6 is used to prove the
following: All left (respectively quasi-) centralizers of $K(A)$ are
locally bounded if and only if $A$ is not isomorphic to $A_{0}\oplus
A_{1}$, for any infinite dimensional elementary $C^{*}$-algebra $A_{1}$.
Also, for any $A$ and any left or quasi-centralizer $T$ of $K(A)$, we
have:
$$
A= A_{0}\oplus A_{1}, \ \text{where} \ A_{1} \ \text{is a dual} \ C^{*} \
\text{-algebra and} \ T \ \text{is
locally bounded on} \ K(A_{0}).\tag1
$$
(The meaning of ``dual" is that $A_{1}$ is a $c_{0}$ direct sum of
elementary $C^{*}$-algebras.)

\noindent{\it Definition 1.6.} If $X$ is a subspace of $A^{**}p$ for an
open projection $p$ and $T$ is a linear transformation from $X$ into
$A^{**}p$, then $T$ is called {\it local} if $xq=0$ implies $(Tx)q=0$ for
all open projections $q$ in $A^{**}$ such that $q\leq p$ and all $x$ in
$X$.

Section 5 deals with $T: X\to A^{**}p\cap A$, where $p$ and $q$ are open
projections and $X=A\cap (qA^{**}p)$; i.e., $X$ is the intersection of a
closed left ideal and a closed right ideal of $A$.  Theorem 5.2 states that
local implies (purely) decomposable in this case.  We then then deduce from
our earlier results:
$$
\aligned
&\text{The ideal generated by} \ X \ \text{can be written as} \ A_{0}\oplus
A_{1}\oplus \dots \oplus A_{n}\\ 
&\text{such that} \ A_{1},\dots, A_{n}  
 \ \text{are elementary and}\\
& \ T_{|A_{0}\cap X} \ \text{is bounded.}\endaligned\tag2
$$

These theorems were inspired by Peetre's  theorem \cite{23, 24}, a special
case of which is:
$$
\aligned
&\text{Any local linear operator on} \ C^{\infty}(X), \ \text{where} \ X \
\text{is a differentiable manifold,}\\
&\text{ is locally given by differential  
operators with smooth coefficients.}\endaligned\tag3
$$
Actually our result is more directly analogous to the following somewhat
easier analogue of Peetre's theorem (in which $X$ is a topological space):
$$
\text{Any local linear operator on} \ C(X) \ \text{is multiplication by a
continuous function.}\tag4
$$

In (3) and (4) $T$ is local if $Tf_{|U}=0$ whenever $U$ is an open subset
of $X$ such that $f_{|U}=0$, or, 
equivalently, if supp $Tf\subset$ supp $f$, where ``supp" denotes closed
support.  There is also a concept of ``decomposable" in the context of
Peetre's theorem:  In the case of (4) decomposable means that $f(x)=0$
implies $(Tf)(x)= 0$ (i.e., delete ``closed" from ``closed support"), and
in the case of (3) it means that $(Tf)(x)=0$ whenever $f$ and its
derivatives of order at most $n_{x}$ vanish at $x$.  Since it is
considerably easier to prove the conclusions of (3) or (4) for
decomposable operators, these results can arguably be interpreted as
statements that local implies decomposable.  Since it is easy to prove that
any local operator which is continuous in a suitable sense is decomposable,
(3) and (4) can also be regarded as automatic continuity theorems.

The special case of (4) where $X= {\Bbb R}$ was a Putnam problem in 1966
\cite{26}, and the general case of (4) is Theorem 9.8 of Luxemburg
\cite{19}.  Theorem 6.3 of Neumann and Ptak \cite{20} generalizes the
locally compact case of (4).  Our technique for proving Theorem 5.2 below
is modelled on a proof of (4) which is different from any of the proofs
cited above (and harder than the proofs in \cite{19} and \cite{26}), but is
partially similar to the proof in \cite{19} of (3) and to the argument on
pages 168-169 of \cite{19} (cf. also Lemma 2 of \cite{23}).  Our earlier
results also use established techniques but are not consciously modelled on
any specific theorem in the automatic continuity literature.

We do not know whether there is a good non-commutative $(C^{*}$-algebraic)
analogue of (3).  Note, though, that the full version of Peetre's
theorem, which deals with transformations from $C^{\infty}$ functions of
compact support to distributions and differential operators with
distribution coefficients, allows a discrete set of exceptional points.  It
seems intriguing to view the dual $C^{*}$-algebras $A_{1}$ of (1) and
$A_{1} \oplus\dots \oplus A_{n}$ of (2) in the same light as this discrete
set, despite the fact that the analogy is not very close.  Of course, any
non-commutative analogue of (3) would deal with Frechet spaces and
algebras rather than Banach spaces and algebras as in the present paper.
The author at one time hoped that such a generalization of this paper would
establish a concept of non-commutative differential operator directly
applicable to unbounded derivations of $C^{*}$-algebras.  Although this
hope now seems wrong, it still appears likely that Section 5 will be
useful.  The reader should compare the concepts of locality and
decomposability used in this paper with the concept of locality defined by
Bratteli, Elliott, and Evans in \cite{6, p.251}.

We are grateful to H. Kim and N-C Wong for stimulating duscussions and to K.
Laursen, M. Neumann, and W. Luxemburg for bibliographical help.
\vfill\eject
\subheading{\S 2. Decomposable linear transformations on certain subspaces
of von Neumann algebras}
\proclaim{Theorem 2.1} Let $M$ be a von Neumann algebra and $p,q$
projections in $M$.  Then the following are equivalent:

\item{(i)}
All decomposable linear transformations from $qMp$ to $M$ are bounded.
\item{(ii)}
If $z$ is a central projection of $M$ such that $zpMp$ is finite
dimensional, then $zqMp$ is finite dimensional.
\endproclaim

{\it Proof.} \ First assume (ii) and let $T: qMp\to M$ be
decomposable.

1) If $p_{1}$ is a subprojection of $p$, then $T$ sends $qMp_{1}$ into
$Mp_{1}$ and $qM(p-p_{1})$ into $M(p-p_{1})$.  Thus $T(xp_{1})= (Tx)p_{1}$
for $x$ in \ $qMp$; i.e., $T$ commutes with $R_{p_{1}}$, the right multiplication
by $p_{1}$.

If $T$ were known to be bounded, it would now be easy to prove that $T$
commutes with all right multiplications by elements of $pMp$.  We proceed
to prove another partial result.

2) If $p_{1}$ and $p_{2}$ are equivalent mutually orthogonal
subprojections of $p$ and $u$ is a partial isometry such that $u^{*}u=
p_{1}$ and $uu^{*} = p_{2}$, then $T$ commutes with $R_{u}$.

To see this, note that $p_{1}, p_{2}, u$, and $u^{*}$ span a copy of
$M_{2}$, the algebra of $2\times 2$ complex matrices, inside $pMp$.  Since
$M_{2}$ is the linear span of its projections, it follows from 1) that $T$
commutes with $R_{x}$ for all $x$ in $M_{2}$.

3) If $\{p_{i}\}$ is a set of mutually orthogonal subprojections of $p$,
then there is a constant $C$ such that for all but finitely many $i,
||T_{|qMp_{i}}||\leq C$.

If this were false, we  could change notation and assume a sequence
$\{p_{n}\}$ such that for each $n, ||T_{|qMp_{n}}|| >n2^{n}$ (this
allows the possibility that $T_{|qMp_{n}}$ is unbounded).  Then choose
$x_{n}$ in $qMp_{n}$ such that $||x_{n}||=2^{-n}$ and $||Tx_{n}|| > n$, and
let $x= \Sigma^{\infty}_{1}x_{n}$.  Since $Tx_{n}= T(xp_{n}) = (Tx)p_{n},
\ ||Tx_{n}||\leq ||Tx||, \ \forall n,$ a contradiction.

4) We now use 1) and 3) to show that there is a largest central projection
$z$ such that $T_{|zqMp}$ is bounded.  By Zorn's lemma there is a maximal
collection $\{z_{i}\}$ of non-zero, mutually orthogonal, central projections
of $M$ such that $T_{|z_{i}qMp}$ is bounded, \ $\forall i$.  

Let $z= \Sigma z_{i}$.  Then $zM$ can be identified with $\underset
i\to\oplus z_{i}M$, the $\ell^{\infty}$ direct sum, and $zqMp$ can be
identified with $\underset i\to\oplus z_{i}qMp$.  Since $T$ commutes with $R_{z_{i}}, \
\forall i, \ T(\underset i\to\oplus x_{i})= \underset i\to\oplus
T_{i}x_{i}$, for $x_{i}$ in $z_{i}qMp$ and $T_{i}= T_{|z_{i}qMp}$.  Since
$\{||T_{i}||\}$  is bounded by 3), it follows that $T_{|zqMp}$ is bounded.
Clearly $z$ is as desired.

5) We now change notation (replace $M$ by $M(1-z))$ and assume $z=0$.  Then
by 3), $M$ does not possess an infinite set  of non-zero, mutually
orthogonal, central projections; i.e., $M$ is the direct sum of finitely
many factors.  Now we easily reduce to the case where $M$ is a factor
and $p,q$ are non-zero.  In  view of (ii), only three cases are possible:

(a) $pMp$  is properly infinite.

(b) $pMp$ is of type II$_{1}$

(c) $qMp$ is finite dimensional.\newline Of course, case (c) is trivial.

6) If $pMp$ is a properly infinite factor, there is a sequence $\{p_{n}\}$
of mutually orthogonal subprojections of $p$, such that each $p_{n}$ is
equivalent to $p$ and $\Sigma_{1}^{\infty}p_{n}=p$.  By (3),
$T_{|qMp_{n}}$ is bounded for at least one $n$.  Since 2) applies to the
pair $p_{n},p-p_{n})$, also $T_{|qM(p-p_{n})}$ is bounded.  Hence $T$ is
bounded.

7) Now assume $pMp$ is a factor of type II$_{1}$ and $T$ is unbounded.
Choose a projection $p_{1}$ in $pMp$ such that $p_{1}$ is equivalent to
$p-p_{1}$.  Applying 2) to the pair $(p_{1},p-p_{1})$, we see that both
$T_{|qMp_{1}}$ and $T_{|qM(p-p_{1})}$ are unbounded .  Next choose a
projection $p_{2}$ in $(p-p_{1})M(p-p_{1})$ such that $p_{2}$ is equivalent
to $p-p_{1}-p_{2}$ and apply the same argument.  Continuing in this manner,
we obtain a sequence $\{p_{n}\}$ of mutually orthogonal subprojections of
$p$ such that $T_{|qMp_{n}}$ is unbounded, \ $\forall n$, contradicting
3).

Now assume (ii) is false.  Then there is central projection $z$ such that
$zM$ is a type I factor, $zpMp$ is finite dimensional and $zqMp$ is
infinite dimensional.  Replacing $M$ with $zM$ and changing notation, we
assume $M= B(H)$ for a Hilbert space $H, \ p$ is a non-zero finite rank
projection on $H$ and $q$ is an infinite rank projection on $H$.  Let $t:
qH\to H$ be a discontinuous linear transformation.  For each $x$ in
$qB(H)p, \ tx$ is a linear transformation whose kernel contains $(1-p)H$.
Therefore $tx\in B(H)$.  Thus we can define an unbounded linear
transformation $T$ by $Tx = tx$ for $x$ in $qMp$, and clearly $T$ is
decomposable.
\proclaim{Corollary 2.2} If $p=q$, then all decomposable linear transformations from $qMp$ to $M$ are bounded.  In particular, all decomposable linear
transformations from $M$ to itself are bounded.
\endproclaim

\noindent {\it Remark 2.3.} It is possible to describe the most general
decomposable linear transformation $T: qMp\to M$.  We have already
mentioned that it follows from results of Wong \cite{27} that if $T$ is
bounded it is given by a left multiplication.  Now by the above proof, the
general case is reduced to the case considered in the last paragraph of the
proof; and using 2), we see that the unbounded decomposable linear
transformations constructed there are the only ones possible.  The
conclusion is that $M,p$, and $q$ can be identified with $M_{0}\oplus
B(H_{1})\oplus\dots\oplus B(H_{n}), \ p_{0}\oplus\dots\oplus p_{n},
q_{0}\oplus\dots\oplus q_{n}$, so that $T$ is the left multiplication
induced by $t_{0}\oplus\dots\oplus t_{n}$, where $t_{0}\in M_{0}q_{0}, \
t_{i}: q_{i}H_{i}\to H_{i}$ for $i > 0$, and $t_{1},\dots, t_{n}$ may be
discontinuous.
 
The  next corollary is due to Wong \cite{27} in the bounded case.
\proclaim{Corollary 2.4} If $M, p$, and $q$ are as above, then a linear
transformation $T:qMp\to Mp$ is decomposable if and only if $T$ commutes
with $R_{x}, \ \forall x\in pMp;$ i.e., $T$ is a homomorphism of right
$pMp$-modules.
\endproclaim
\vfill\eject
\subheading{\S 3. Purely decomposable linear transformations in the context
of $C^{*}$-algebras}

The decomposable linear transformations which occur in Wong \cite{27} arise
in connection with a closed projection $p$, but the analysis of centralizers
of Pedersen's ideal instead uses open projections.  In order to cover both
cases, our basic theorem is stated in an abstract way.  If $A$ is a
$C^{*}$-algebra and $\pi: A\to B(H)$ a (*-) representation, $\pi$ can be
uniquely extended to a normal representation of $A^{**}$, the enveloping
von Neumann algebra.  This extension will be denoted by $\pi^{**}$.

\noindent {\it Definition 3.1.} If $A$ is a $C^{*}$-algebra, $p_{1}$ and
$p_{2}$ are projections in $A^{**}$, and $X$ is a subspace of
$p_{1}A^{**}p_{2}$, then $X$ is said to satisfy $K(p_{1},p_{2})$ if the
following is true: If $\pi: A\to B(H)$ is any irreducible representation
such that $\pi^{**}(p_{1}), \pi^{**}(p_{2})\not= 0$, and if $V$ is any
non-trivial finite dimensional subspace of $\pi^{**}(p_{2})H$, then there
is a positive number $k(\pi,V)$ such that for any $t$ in
$\pi^{**}(p_{1})B(H)$, there is an $x$ in $X$ with $||x||\leq
k(\pi,V)||t||$ and $\pi^{**}(x)_{|V} = t_{|V}$.  If $X$ satisfies
$K(p_{1},p_{2})$ and $\pi$ is as above, then $k(\pi)$ denotes
inf$\{k(\pi,V): V$ is a 1-dimensional subspace of $\pi^{**}(p_{2})H.\}$

\noindent {\it Definition 3.2.} If $p_{1}$ and $p_{2}$ are each open or
closed, a subspace $X(p_{1},p_{2})$ is defined as follows:

\item{(a)}
If $p_{1}$ and $p_{2}$ are closed, $X(p_{1},p_{2}) = p_{1}Ap_{2}$.

\item{(b)} 
If $p_{1}$ and $p_{2}$ are open, $X(p_{1},p_{2})= A\cap
(p_{1}A^{**}p_{2})$.

\item{(c)}
If $p_{1}$ is open and $p_{2}$ is closed,$X(p_{1},p_{2}) = Ap_{2}\cap
(p_{1}A^{**})$.

\item{(d)}
If $p_{1}$ is closed  and $p_{2}$ is open, $X(p_{1},p_{2})= p_{1}A\cap
(A^{**}p_{2})$.

Then $X(p_{1},p_{2})$ is always norm closed.  See 4.4 of \cite{5} (cf. also
pages 916-918 of \cite{8}) for how to deal with closed projections.  In
case (a), $p_{1}Ap_{2}$ is isometric to $A/L+R$, where $L$ is the left
ideal corresponding to $p_{2}, R$ the right ideal corresponding to $p_{1}$,
and $L+R$ is closed by a result of Combes \cite{11}.  If $p_{1}$ and/or
$p_{2}$ is both open and closed, it is not hard to show that both
definitions of $X(p_{1},p_{2})$ agree.  (A projection is both open and
closed if and only if it is a multiplier of $A$.)
\proclaim{Proposition 3.3} If $X(p_{1},p_{2})$ is as defined above, then
$X(p_{1},p_{2})$ satisfies $K(p_{1},p_{2})$ with $k(\pi,V)=1, \ \forall \
\pi,V$.
\endproclaim

{\it Proof.} Let $\pi, V$, and $t$ be as in 3.1.  By the Kadison
transitivity theorem \cite{16}, there is $a$ in $A$ such thate $||a||\leq
||t||$ and $\pi(a)_{|V} = t_{|V}$.  Let $q_{2}$ be the projection with
range $V$ and $q_{1}$ the projection with range $t(V)$.  If $p_{i}$, for
$i=1$ or $2$, is open then Akemann's Urysohn lemma \cite{3, Lemma III.1}
yields $b_{i}$ in $A$ such that $b_{i}= p_{i}b_{i}p_{i}, q_{i} = \pi(b_{i})q_{i},$ and $||b_{i}||=1$.  If $p_{i}$ is closed, let $b_{i}= p_{i}$.  Then in all
cases let $x = b_{1}ab_{2}$, and it is easy to see that $x\in
X(p_{1},p_{2}), ||x||\leq ||a||\leq ||t||$, and $\pi^{**}(x)_{|V}=
\pi(a)_{|V}= t_{|V}$.

\noindent {\it Remarks 3.4.} (a). In case (b) of 3.2, $X(p_{1},p_{2})$ is a Hilbert $B_{1}-B_{2}$ bimodule, where $B_{i}$ is the hereditary $C^{*}$-subalgebra of $A$ supported by $p_{i}$ (notation: $B_{i} =  her(p_{i}))$.  Also $X(p_{1},p_{2})= B_{1}AB_{2}= (B_{1}AB_{2})^{-}$, and $X(p_{1},p_{2}) = R\cap L$, where $R= A\cap (p_{1}A^{**})$, the right ideal corresponding to
the open projection $p_{1}$ (which is the same as the right ideal
corresponding to the closed projection $1-p_{1}$), and $L$ is the left
ideal corresponding to the open projection  $p_{2}$.  Also
$her(p_{i})=X(p_{i},p_{i}), L= X(1,p_{2})$, and $R= X(p_{1},1)$. 

(b) In case (c) of 3.2 it can be shown that $X(p_{1},p_{2})= Rp_{2}$ where
$R$  is as above; and a similar alternate description can be given in case
(d) of 3.2.  We have no application in mind for these mixed cases but have
included them for completeness, with the feeling that they could prove
useful.

(c) The main hypothesis on the domain of the decomposable transformation in
the next theorem is $K(p_{1},p_{2})$, and the exact form of Definition 3.1
represents a compromise on the part of the author.  On the one hand
$K(p_{1},p_{2})$ is not a minimal hypothesis, but it is more simply stated
than a minimal hypothesis and the additional generality gained from using a
weaker hypothesis is probably not important.  Also, if $X$ is norm closed,
the existence of $k(\pi,V)$ could be deduced from the open mapping
theorem.  On the other hand Definition
3.1 could be simplified by requiring $k(\pi,V)=1$, and the generality lost
by thus strengthening the hypothesis is probably not important.

(d) In contrast to the above, the focus on irreducible representations in
this section is not an arbitrary decision of the author.  It seems
necessary to use the Kadison transitivity theorem to obtain an applilcable
automatic continuity theorem in this context, and this forces us to use
irreducible representations.

There is one more technical point needed.  Let $\pi_{a}$, the reduced atomic
representation of $A$, be the direct sum of one irreducible
representation from each unitary equivalence class.  Thus for $x$ in
$A^{**}, \ ||\pi^{**}_{a}(x)|| = \sup\{||\pi^{**}(x)||: \pi$ irreducible\}.
Pedersen \cite{21} defined a large, norm-closed space $\Cal U$, consisting of
self-adjoint elements of $A^{**}$ and called the set of universally
measurable operators, and showed that $||\pi^{**}_{a}(x)||= ||x||, \ \forall
x\in\Cal U$ (cf.  \cite{22, \S 4.3}).  Because variants of this concept were
needed in Wong \cite{27}, it seems desirable to abstract this property
also.

\noindent {\it Definition 3.5.} If $Y$ is a subspace of $A^{**}, Y$ is of
{\it type} $\Cal U$ if $\pi^{**}_{a}(y)=0$ implies $y=0$ for all $y$ in $\overline
Y$, the
norm closure of $Y$.  If there is a constant $c$ such that $||y||\leq
c||\pi^{**}_{a}(y)||$, for all $y$ in $Y$ (or, equivalently, in $\overline
Y), Y$ is of {\it strong type} $\Cal U$.

Thus any subspace of $\Cal U_{{\Bbb C}}$, the complexification of $\Cal
U$, is of
strong type $\Cal U$.  Unfortunately it is not known whether $\Cal U_{{\Bbb C}}$ is
an algebra.  A future paper of the author \cite{9} will show that 
$\Cal U$ is the $C^{*}$-algebra generated by $A$ and all open or closed
projections.  Wong \cite{27}  uses some variants of the concept of
universal measurability and shows the appropriate spaces are of type
$\Cal U$.  There are a number of results in the literature to the
effect that $\Cal U$ or $\Cal U_{{\Bbb C}}$ is ``sufficiently large" (these include
\cite{12, 2.2.15, 2.4.3}, \cite{21, Prop. 3.6}, \cite{7, Thm. 4.15},
\cite{8. Thm. 4.10}, and \cite{4, Lemma 2.1}).  Also, if $A$ is separable,
$\{x\in A^{**}: x$ satisfies the barycenter formula\} is a weakly sequentially	
closed $C^{*}$-algebra of strong type  $\Cal U$ which contains 
$\Cal U_{{\Bbb C}}$.
This algebra is large enough for all imaginable applications.  The upshot
is that the ``type $\Cal U$" hypothesis in the next theorem does not impede
applicability.  
\proclaim{Theorem 3.6} Assume $A$ is a $C^{*}$-algebra, $X$ is a subspace
of $A^{**}$ which satisfies $K(p_{1},p_{2})$, where $p_{1}$ and $p_{2}$
are projections in $A^{**}$, and $T: X\to A^{**}$ is a purely decomposable
linear transformation.

(a) For each irreducible representation $\pi: A\to B(H_{\pi})$ such that
$\pi^{**}(p_{1}), \pi^{**}(p_{2})\not= 0$, there is a linear transformation
$t_{\pi}: \pi^{**}(p_{1}) H_{\pi}\to H_{\pi}$ such that $\pi^{**}(Tx) =
t_{\pi}\pi^{**}(x), \ \forall \ x\in X$.  Also if $\pi^{**}(p_{1})$ or
$\pi^{**}(p_{2})$ is $0$, then $\pi^{**}[T(X)]= \{0\}$.

(b) If $T$ is bounded, then $t_{\pi}$  is bounded and $||t_{\pi}||\leq
k(\pi) ||T||$.

(c) If $X$ is norm closed and either $\pi^{**}(p_{2})$ has infinite rank or
$\pi^{**}(p_{1})$ has finite rank, then $t_{\pi}$ is bounded.

(d) Assume that $X$ is norm closed, $T(X)$ is of type $\Cal U$, and for each
irreducible representation $\pi$ of $A$ at least one of the following
holds:

(i) $\pi^{**}(p_{2})=0$

(ii) $\pi^{**}(p_{2})$ has infinite rank

(iii) $\pi^{**}(p_{1})$ has finite rank.

\noindent Then $T$ is bounded.
\endproclaim
\noindent {\it Remark.} The main part of the sufficient condition for
automatic continuity given by part (d) is analogous to that in Theorem 2.1,
but it is not necessary.  In the following two sections, in the context of
open projections $p_{1}, p_{2}$, we will obtain necessary and sufficient
conditions for automatic continuity, but we have not attempted this in the
context of closed projections.

{\it Proof.} (a) If $v\in H_{\pi}, \ x\in X$, and $\pi^{**}(x)v = 0$, then
$\pi^{**}(Tx)v=0$, since $T$ is purely decomposable.
($||\pi^{**}(a)v||^{2}= \varphi_{v}(a^{*}a), \ \forall a\in A^{**}$, where
$\varphi_{v}$ is a multiple of a pure state, $\varphi_{v}=
(\pi(\cdot)v,v).)$.  Thus we can define a linear transformation $t_{v}$
on $\{\pi^{**}(x)v: x\in X\}$ by $t_{v}\pi^{**}(x)v= \pi^{**}(Tx)v$.  If
$v\notin \pi^{**}(1-p_{2})H_{\pi}$, it follows from $K(p_{1},p_{2})$ that
the domain of $t_{v}$ is $\pi^{**}(p_{1})H_{\pi}$
 and clearly $t_{cv+w}= t_{v}$ for
non-zero $c$ in ${\Bbb C}$ and $w$ in $\pi^{**}(1-p_{2})H_{\pi}$.  If
$v_{1}$ and $v_{2}$ are linearly independent modulo
$\pi^{**}(1-p_{2})H_{\pi}$ and $w$ is in $\pi^{**}(p_{1})H_{\pi}$, then by
$K(p_{1},p_{2})$ there is $x$ in $X$ such that $\pi^{**}(x)v_{1}=
\pi^{**}(x)v_{2}=w$.  Thus $\pi^{**}(x)(v_{1}-v_{2})=0$ and since $T$ is
purely decomposable, $\pi^{**}(Tx)(v_{1}-v_{2})=0$.  Since $t_{v_{i}}w=
\pi^{**}(Tx)v_{i}$, this implies $t_{v_{1}}= t_{v_{2}}$.  Thus we can take
$t_{\pi}$ to be the common value of $t_{v}$.  The last sentence is clear.

(b) Let $v$ be a unit vector in $\pi^{**}(p_{2})H_{\pi}$ and $V= {\Bbb
C}v$.  For any $w$ in $\pi^{**}(p_{1})H_{\pi}$, there is $x$ in $X$ such
that $\pi^{**}(x)v=w$ and $||x||\leq k(\pi,V)||w||$.  Then
$||t_{\pi}w||=||\pi^{**}(Tx)v|| \leq ||T|| ||x||\leq k(\pi,V)||T|| ||w||$.
Thus $||t_{\pi}||\leq k(\pi,V)||T||$ for all such $V$, and hence
$||t_{\pi}||\leq k(\pi)||T||$.

(c) It is trivial that $t_{\pi}$ is bounded if $\pi^{**}(p_{1})$ has
finite rank.  Thus assume $\pi^{**}(p_{2})$ has infinite rank and
$t_{\pi}$ is unbounded.  Choose an orthonormal sequence $\{e_{n}\}$ in
$\pi^{**}(p_{2})H_{\pi}$, and let $V_{n} = span \{e_{1},\dots, e_{n}\}$.
We choose recursively $x_{n}$ in $X$ such that $||x_{n}|| < 2^{-n},
\pi^{**}(x_{n})e_{k}=0$ for $k < n$,  and $||t_{\pi}\pi^{**}(x_{1}+ \dots +
x_{n})e_{n}|| > n$.  If $x_{k}$ has been chosen for $k < n$, we choose $w$
in $\pi^{**}(p_{2})H_{\pi}$ such that $||w|| < 2^{-n}k(\pi,V_{n})^{-1}$ and
$||t_{\pi}w|| > n + ||t_{\pi}\pi^{**}(x_{1}+\dots x_{n-1})e_{n}||$.  Then by
$K(p_{1},p_{2})$ we can find $x_{n}$ in $X$ such that
$\pi^{**}(x_{n})e_{n}=w, \pi^{**}(x_{n})e_{k}=0$ for $k<n$, and
$||x_{n}||\leq k(\pi,V_{n})||w||< 2^{-n}$.  Now let $x=
\Sigma^{\infty}_{1}x_{k} = y_{n} + z_{n},$ where $y_{n}=
\Sigma^{n}_{1}x_{k}$ and $z_{n} = \Sigma^{\infty}_{n+1}x_{k}$.  Then
$\pi^{**}(z_{n})e_{n}=0$.  Hence $||t_{\pi}\pi^{**}(x)e_{n}||=
||t_{\pi}\pi^{**}(y_{n})e_{n}|| > n$.  Then $t_{\pi}\pi^{**}(x)$ is
unbounded, which contradicts the fact that $t_{\pi}\pi^{**}(x) =
\pi^{**}(Tx)$.

(d) This now follows from the closed graph theorem.  If $T$ is not
bounded, there is a sequence $\{x_{n}\}$ in $X$ such that $x_{n}\to 0$ and
$Tx_{n}\to y$, for some non-zero $y$ in $T(X)^{-}$.  Then there is an
irreducible $\pi$ such that $\pi^{**}(y)\not= 0$.  By the last sentence of
(a) this implies $\pi^{**}(p_{1}), \pi^{**}(p_{2})\not= 0$.  Hence $t_{\pi}$
is defined and by (c) it is bounded.  Hence $\pi^{**}(y)=
\lim\pi^{**}(Tx_{n})= \lim t_{\pi}\pi^{**}(x_{n})=0$, a contradiction.

Wong \cite{27} considers a closed projection $p$ in $A^{**}$ and two spaces
of ``sections", denoted $\Cal S$ and $\Cal W$.  $\Cal S$ can be
identified with $Ap$ and $\Cal W$ with $\{x\in A^{**}p: y^{*}x\in pAp,
\ \forall \ y\in Ap\}$.
\proclaim{Corollary 3.7} If there is no infinite dimensional irreducible
representation $\pi$ of $A$ such that $\pi^{**}(p)$ has non-zero finite
rank, then any purely decomposable linear transformation from $\Cal S$
to $\Cal W$ or from $\Cal W$ to $\Cal W$ is bounded.  ({\it A
fortiori}, the same is true for transformations from $\Cal S$ or
$\Cal W$ to $\Cal S$.)
\endproclaim

{\it Proof.} Since $\Cal S= X(1,p)$ and $\Cal W\supset \Cal S,
\Cal S$ and $\Cal W$ satisfy $K(1,p)$ by 3.3, and both are norm
closed.  $\Cal W$ is of type $\Cal U$ by \cite{27}.
\proclaim{Corollary 3.8} If $p$ is a closed projection in 
$A^{**}$ and $T: pAp\to \Cal U_{\Bbb C}$ is a purely decomposable linear
transformation, then $T$ is bounded.  In partilcular, any purely
decomposably linear transformation from $A$ to $\Cal U_{{\Bbb C}}$ is
bounded.
\endproclaim

\noindent{\it Remark 3.9.} If $T$ is decomposable, $T(X)\subset
A^{**}p_{2}$.  We briefly discuss the question of proving {\it a priori}
that purely decomposable implies decomposable.  If $x$ is in $X$ and $q=
\In \ x$, pure decomposability implies that $\pi^{**}_{a}[(Tx)(1-q)]=0$.  If
$(Tx)(1-q)$ is contained in some space of type $\Cal U$, then this
implies that $(Tx)(1-q)=0$; i.e. $\In \ Tx\leq \In \ x$.  Now if the $C^{*}$-algebra generated by $x$ is contained in $\Cal U_{{\Bbb C}}$ then
$q\in\Cal U_{{\Bbb C}}$ by \cite{22, 4.5.15}.  Under suitable universal
measurability hypotheses on the elements of $X$ and $T(X)$, this sort of
argument can be used to prove decomposability (cf. the comments after
Definition 3.5).  If we want only to show that $T(X)\subset A^{**}p_{2}$,
it is enough to assume (or prove) $T(X)(1-p_{2})$ is of type $\Cal U$.

It has already been mentioned that in many cases a bounded purely
decomposable linear transformation can be proved to be a left
multiplication.  This in particular would be an {\it a postiori} proof of
decomposability.  It is proved in \cite{27} that transformations of the
four types mentioned in 3.7 are left multiplications (if bounded).  We
content ourselves here with proving the easiest abstract result of this
sort.  Of course this is not an automatic continuity result, but when
applicable it allows the conclusion of Theorem 3.6(d) to be strengthened.
\proclaim{Proposition 3.10} Assume $A$ is a $C^{*}$-algebra, $X$ is a
subspace of $A^{**}$ which satisfies $K(p_{1},p_{2})$ with $k(\pi)$
bounded independently of $\pi$, where $p_{1}$ and $p_{2}$ are projections
in $A^{**}, T: X\to A^{**}$ is a bounded, purely decomposable linear
transformation, and $span(\{(Tx_{1})x_{2} \dots x_{n}: x_{i}\in X$ for $i$
odd, $x_{i}\in X^{*}$ for $i$ even, $n\geq 1\})$ is of strong type $\Cal U$.  Then
there is $t$ in $A^{**}$ such that $Tx = tx, \ \forall \ x\in X$.
\endproclaim

{\it Proof.} Let $Y = span(\{x_{1}x_{2}\dots x_{n}: x_{i}\in X$ for $i$ odd,
$x_{i}\in X^{*}$ for $i$ even, $n\geq 1\}$) and let $Z$ be the linear span
above.  Define a linear map $S: Y\to Z$ by $S(x_{1}x_{2}\dots x_{n})=
(Tx_{1})x_{2}\dots x_{n}$.   We show that $S$ is bounded, and hence
well-defined.  For $\pi$ irreducible with $\pi^{**}(p_{1}),
\pi^{**}(p_{2})\not= 0$ and $y \ \in \ Y, ||\pi^{**}(Sy)||=
||t_{\pi}\pi^{**}(y)||\leq ||t_{\pi}||\cdot||y||\leq k(\pi)||T||\cdot ||y||\leq
c_{1}||T||\cdot ||y||,$ by 3.6 (a), (b) and hypothesis.  If
$\pi^{**}(p_{1})$ or $\pi^{**}(p_{2})$ is $0$, then $\pi^{**}(X)=
\pi^{**}(T(X))= \{0\},$ and hence $\pi^{**}(Sy)=0$.  Thus $||Sy||\leq
c_{2}||\pi^{**}_{a}(Sy)||\leq c_{1}c_{2}||T|| \ ||y||$, since $Z$ is of
strong type $\Cal U$.   Then $S$ extends to a bounded linear map from
$\overline Y$ to $\overline Z$, also denoted by $S$.  Let $\{e_{i}\}$ be an
approximate identity of the $C^{*}$-algebra generated by $XX^{*}$.  Then
for $x$ in $X, \ e_{i}x\in \overline Y, \ S(e_{i}x)= S(e_{i})x,$ and
$e_{i}x\to x$ in norm.  Hence $Tx= Sx = \lim S(e_{i}x)= \lim S(e_{i})x$.
Now since $||Se_{i}||\leq c_{1}c_{2}||T||, \ \{Se_{i}\}$ has a weak cluster
point $t$ in the von Neumann algebra $A^{**}$.  Then $Tx = tx, \ \forall \
x\in X$, and also $||t||\leq c_{1}c_{2}||T||$. (The earlier results and
remarks of this section show that in many cases $c_{1}c_{2}=1$.)
\vfill\eject
\subheading{\S 4. Centralizers of Pedersen's ideal}

If $B$ is an algebra, a {\it left centralizer} of $B$ is a linear
transformation $T: B\to B$ such that $TR_{b}=  R_{b}T$ for every right
multiplication $R_{b}, \ b$ in $B$.  Right centralizers are defined
similarly.  A {\it quasi-centralizer} is a bilinear map $T: B\times B\to B$
such that $T(b_{1}b_{2}, b_{3})= b_{1}T(b_{2},b_{3})$ and
$T(b_{1},b_{2}b_{3}) = T(b_{1},b_{2})b_{3}$ for all $b_{1},b_{2}, b_{3}$
in $B$.  A {\it double cenralizer} is a pair $(T_{1},T_{2})$ of linear maps
such that $b_{1}(T_{1}b_{2})= (T_{2}b_{1})b_{2}$ for all $b_{1}, b_{2}$ in
$B$.  If the left and right annihilators of $B$ are trivial, it then
follows that $T_{1}$ is a left centralizer and $T_{2}$ a right
centralizer.  Note that $B$ should be assumed non-unital, since the theory
of centralizers is trivial if $B$ is unital.

If $C$ is a super-algebra of $B$, then a {\it left multiplier} of $B$ (in
$C$) is an element $c$ of $C$  such that $cB\subset B$.  Similarly $c$ is a
{\it right multiplier} if $Bc\subset B$, a {\it quasi-multiplier} if
$BcB\subset B$, and a {\it multiplier} if  $cB, Bc\subset B$.  Every left
multiplier induces a  left centralizer (namely, $L_{c|B})$, every
quasi-multiplier induces a quasi-centralizer, etc.

If $A$ is a $C^{*}$-algebra, when $B$ above is taken to be $A$, it is
standard to take $C$ to be $A^{**}$.  In this case the sets of multipliers
of the various types are denoted $LM(A), \ RM(A), \ QM(A)$, and $M(A)$.  Of
course, $M(A)= LM(A)\cap RM(A)$ and $LM(A), \ RM(A)\subset QM(A)$.  Also,
every left centralizer of $A$ is induced from a left multiplier, every
quasi-centralizer from a quasi-multiplier, etc., and moreover the norms
agree.  If $\pi: A\to B(H)$ is any faithful, non-degenerate
representation, then $\pi^{**}$ maps $LM(A)$ isometrically onto the set of
left multipliers of $\pi(A)$ in $B(H)$, and similarly  for $RM(A), QM(A),
M(A)$.  References for the above, which has been sketched only very
briefly, are Johnson \cite{15}, Busby \cite{10}, and section 3.12 of \cite{22}.

In this section the role of $B$ above is played by $K(A)$, Pedersen's
ideal, and centralizers of the various types may not be bounded.  Of
course, any bounded centralizer extends by continuity to a centralizer of
$A$.  If $A$ is the commutative $C^{*}$-algebra $C_{0}(X)$, then the
centralizers (of any type) of $A$ can be identified with $C_{b}(X)$, the
set of bounded continuous functions on $X$, and the centralizers (of any
type) of $K(A)$ can be identified with $C(X)$, the set of arbitrary
continuous functions.  

Lazar and Taylor \cite{18} and Phillips \cite{25} showed that any double
centralizer $T= (T_{1},T_{2})$ of $K(A)$ is locally bounded in the
following sense: There is a family $\{I_{j}\}$ of (closed two-sided)
ideals of $A$ such that $A= (\Sigma I_{j})^{-}$ and $T_{|K(A)\cap I_{j}}$
is bounded for each $j$.  For each $j, T$ induces an element of $M(I_{j})$;
and $M(I_{j})$ will be identified with a subset of $A^{**}$, since
$A^{**}$ is canonically isomorphic to $I^{**}_{j}\oplus (A/I_{j})^{**}$.
Phillips \cite{25} identified $\Gamma(K(A))$, the set of double
centralizers of $K(A)$, with the inverse limit of $\{M(I_{c})\}$ for a
suitable directed family $\{I_{c}\}$ of ideals of $A$.  It is also possible
to identify $\Gamma(K(A))$ with an appropriate set of unbounded operators
(on the Hilbert space of the universal representation of $A$) affiliated
with the von Neumann algebra $A^{**}$, and this approach to $\Gamma(K(A))$
and related concepts is taken by H. Kim \cite{17}.

Kim asked whether quasi-centralizers of $K(A)$ would also be locally
bounded and then answered this question negatively by showing that it is
false for infinite dimensional elementary $C^{*}$-algebras.  (See 4.6 below
and note that since elementary $C^{*}$-algebras are simple, local
boundedness would imply boundedness.)  Infinite dimensional elementary
$C^{*}$-algebras are the easiest examples of non-unital, non-commutative
$C^{*}$-algebras, but they also turn out to be essentially the only
counter-examples for this question.

In this section we will frequently use some widely known facts about
$C^{*}$-algebras concerning the spectrum of a $C^{*}$-algebra $A$ (the
spectrum, denoted by $\hat A$, is the set of equivalence clases of
irreducible representations) and its topology and $CCR$ algebras  (also
called liminal $C^{*}$-algebras).  References for this are sections 3.1 to
3.5, 4.2, and 4.4 of \cite{13} or sections 4.1, 6.1, and 6.2 of \cite{22}.
This material, as well as material on the Pedersen ideal contained in
section 5.6 of \cite{22} may be used without explicit citation.
\proclaim{Proposition 4.1} $K(A)$, for a $C^{*}$-algebra $A$, satisfies
$K(1,1)$ with $k(\pi,V)=1, \ \forall \ \pi,V$.
\endproclaim

{\it Proof.} Let $\pi: A\to B(H)$ be irreducible, let $V$ be a
finite-dimensional subspace of $H$, and let $t$ be in $B(H)$.  By the
Kadison transitivity theorem there is $a$ in $A$ such that $||a||\leq
||t||$ and $\pi(a)_{|V} =t_{|V}$.  By the same theorem there is $b$ in $A$
such that $b^{*}=b, \pi(b)v=v, \ \forall \ v\in V$, and $||b|| = 1$.  Let
$f: {\Bbb R}\to [0,1]$ be a continuous function such that $f(1)=1$ and $f$
vanishes in a neighborhood of $0$.  If $c=f(b)$, then $c\in K(A)$ and $c$
satisfies the properties given above for $b$.  Then for $x= ac, \ x\in
K(A), ||x||\leq ||a||\leq ||t||$, and $\pi(x)_{|V}=t_{|V}$.
\proclaim{Proposition 4.2} If $T$ is a left centralizer of $K(A)$, then $T$
is decomposable and {\it a fortiori} purely decomposable.
\endproclaim

\noindent {\it Remark.} It is also true that any purely decomposable $T:
K(A)\to K(A)$ is a left centralizer.  The proof is similar to that of
$(iv)\Rightarrow (i)$ in 5.5 below.

{\it Proof.} If $x$ is in $K(A)$, then $L\subset K(A)$, for $L$ the closed
left ideal generated by $x$, by \cite{18, Prop. 3.3}.  There is $y$ in $L$
such that $x=y(x^{*}x)^{\frac{1}{4}}$ by \cite{22. 1.4.5}.  Since also
$(x^{*}x)^{\frac{1}{4}}\in K(A), \ Tx= (Ty)(x^{*}x)^{\frac{1}{4}}$.
Therefore $\In \ Tx\leq \In(x^{*}x)^{\frac{1}{4}} = \In \ x$.
\proclaim{Lemma 4.3} Assume $L$ is a closed left ideal of the
$C^{*}$-algebra $A, L$ generates $A$ as a (closed, two-sided) ideal, and
$T: L\to L$ is a purely decomposable linear map.  If $T$ is not bounded,
then there is a direct sum decomposition, $A=A_{0}
\oplus A_{1},$ such that $A_{1}$ is an elementary $C^{*}$-algebra for which
the corresponding $t_{\pi}$ (notation as in 3.6) is not bounded.
\endproclaim

\noindent{\it Remark.} In the next section lemmas 4.3 and 4.4 will be
generalized by replacing $L$ with $X(p_{1},p_{2})$ for open projections
$p_{1}$ and $p_{2}$.  The same arguments apply except that some details
must be added to the ``preliminary" paragraph in the proof below dealing
with cutting down to an ideal.  These additional details will be provided
in the next section.

{\it Proof.} Some explanation of the statement of the lemma is in order.
First, $L= X(1,p)$ for some open projection $p$ and hence $t_{\pi}$ is
defined for each irreducible $\pi$.   (Note: $\pi^{**}(p)=0$ implies
$\pi(L)=0$, which is impossible by hypothesis.)  Second, if $A= A_{0}\oplus
\Cal K(H)$, then an irreducible representation $\pi: A\to B(H)$ is
defined as projection onto the second summand.  Obviously, $\pi(A)= \Cal
K(H)$.

Now because of the one-to-one correspondence between ideals of $A$ and
open subsets of $\hat A$, direct sum decompositions of $A$ correspond
one-to-one to separations of $\hat A, \hat A = S_{0}\cup S_{1}$.  If
$A_{1}$, where $\hat A_{1}=S_{1}$, is to be elementary, then $S_{1}$ must
consist of a single clopen point.  Conversely, if $\pi$ is an irreducible
representation such that $\{[\pi]\}$, where $[\pi]$ denotes the
equivalence class of $\pi$, is a clopen subset of $\hat A$ and $\pi(A)=
\Cal K(H_{\pi})$, then $\pi$ corresponds to a direct sum decomposition of
the desired type.  (The last condition is necessary because it is not known
whether  $\hat A_{1}=$ \{one point\} implies $A_{1}$  elementary.)  

In the course of the proof we will replace $A$ by an ideal, and therefore
we include one more preliminary paragraph.  If $I$ is an ideal of $A$, then
$L\cap I$ is a closed left ideal of $I$ and $L\cap I$ generates $I$ as an
ideal (for example, because $I= IAI = I(ALA)^{-}I\subset
(IALAI)^{-}\subset [(IL)I]^{-}$ and $IL\subset L\cap I)$.  Also, $T$ maps
$L\cap I$ into itself, since $L\cap I= L\cap [\cap\{kernel \ \pi: [\pi]\in
(A/I)^{\wedge}\}]$ and $T$ maps $L\cap kernel \ \pi$ into itself by pure
decomposability.  Now for $\pi$ in $\hat I, \ t_{\pi}$ is the  same whether
computed relative to $L\cap I$ or $L$ (the domain of $t_{\pi}$ is
$H_{\pi}$).  Finally, any direct sum decomposition, $I= I_{0}\oplus\Cal
K(H_{\pi})$, leads also to $A=  A_{0}\oplus\Cal K(H_{\pi})$ provided
$\pi(A)= \Cal K(H_{\pi})$.  The reason is that $\pi(A)= \Cal
K(H_{\pi})$ implies $[\pi]$ is closed in $\hat A$, and $[\pi]$ open in
$\hat I$ implies $[\pi]$ open in $\hat A$.

Now since $T$ is not bounded, by the closed graph theorem there is a
sequence $\{x_{n}\}$ in $L$ such that $x_{n}\to 0$ and $Tx_{n}\to y\not=
0$.  Of course, $y$ is in $L$, and since $\pi(y) = \lim\pi(Tx_{n}) = \lim
t_{\pi}\pi(x_{n}), \ t_{\pi}$ must be unbounded for all $\pi$ such that
$\pi(y)\not= 0$.  If $I$ is the ideal generated by $y$, this means that
$t_{\pi}$ is unbounded whenever $[\pi]$ is in $\hat I$.  By 3.6 (c), if
$t_{\pi}$ is not bounded, then $\pi^{**}(p)$ has finite rank and hence
$\pi(x)$ has finite rank for all $x$ in $L$.  Since $L$ generates $A$ as an
ideal, this implies $\pi(A)\subset\Cal K(H_{\pi}),$ which implies
$\pi(A)= \Cal K(H_{\pi})$.

We replace $A$ by $I$ and change notation.  Thus we now have that $t_{\pi}$
is unbounded for all $\pi$ and that $A$ is a CCR algebra.  Then $A$ has a
non-zero ideal which has Hausdorff spectrum, and we may cut down and change
notation once more.  Thus, now $\hat A$ is Hausdorff and $t_{\pi}$ is
unbounded for all $\pi$.

To complete the proof, we need only show that $\hat A$ has an isolated
point, and we do this by showing that $\hat A$ is finite.  In fact, if
$\hat A$ has infinitely many points, then there is a sequence $\{U_{n}\}$
of non-empty, mutually disjoint open sets.  If $I_{n}$ is the ideal
corresponding to $U_{n}$, then $t_{\pi}$ is unbounded for all irreducible
representations $\pi$ of $I_{n}$ (i.e., all irreducibles $\pi$ of $A$ such
that $[\pi]\in U_{n}$, or equivalently $\pi(I_{n})\not= \{0\})$.  By 3.6
(b) this means $T_{|L\cap I_{n}}$ is not bounded.  Thus we may choose
$x_{n}$ in $L\cap I_{n}$ such that $||x_{n}|| < 2^{-n}$ and $||Tx_{n}|| >
n$.  Let $x= \Sigma^{\infty}_{1}x_{n}$, an element of $L$.  For $[\pi]$ in
$U_{n}, \pi(Tx)= t_{\pi}\pi(x)= t_{\pi}\pi(x_{n})$, since $\pi_{|I_{m}}=0$
for $m\not= n$.  Thus $||Tx||\geq \underset{U_{n}}\to\sup
||t_{\pi}\pi(x_{n})||= \underset{U_{n}}\to\sup ||\pi(Tx_{n})||= ||Tx_{n}||
> n,$ where one of the equalities uses the fact that $Tx_{n}\in I_{n}$.
This contradiction completes the proof.
\proclaim{Corollary 4.4} Assume $L$ is a closed left ideal of the
$C^{*}$-algebra $A, L$ generates $A$ as an ideal, and  $T: L\to L$ is a
purely decomposable linear map.  Then there is a finite set $F=
\{[\pi_{1}],\dots, [\pi_{n}]\}$ of clopen points of $\hat A$ such that
$\pi_{i}(A)= \Cal K(H_{\pi_{i}})$ and $t_{\pi}$ is bounded whenever
$[\pi]\notin F$. 
\endproclaim

{\it Proof.}  Let $U= \{[\pi]: [\pi]$ \ is a clopen point of
$\hat A,
\pi(A) = \Cal K(H_{\pi}),$ \ and \ $t_{\pi}$ \ is not
bounded\}.  If 
$U$ is infinite, there is a sequence $\{[\pi_{n}]\}$ of
distinct points in $U$.  Then the argument in the last paragraph of the
proof of 4.3 produces a contradiction, with $U_{n} = \{[\pi_{n}]\}$ and
$I_{n}$ the elementary direct summand corresponding to $\pi_{n}$.  Thus $U$
is finite, and we have a direct sum decomposition,  $A= A_{0}\oplus \Cal
K(H_{\pi_{1}})\oplus\dots \oplus \Cal K(H_{\pi_{n}})$, where $U= \{[\pi_{1}],\dots, [\pi_{n}]\}$.  The proof of
4.3 shows that the hypotheses are still satisfied if we replace $(A,L)$ by
$(A_{0}, A_{0}\cap L)$; and since $U\cap \hat A_{0} = \emptyset$, 4.3 then
implies that $T_{|A_{0}\cap L}$ is bounded.  Then 3.6 (b) implies that
$t_{\pi}$  is bounded for $\pi \not\cong \pi_{1},\dots, \pi_{n}$.
\proclaim{Theorem 4.5} If $A$ is a $C^{*}$-algebra and $T$ is a left
centralizer of $K(A)$, then there is a direct sum decomposition, $A=
A_{0}\oplus A_{1}$,  such that $A_{1}$ is a dual $C^{*}$-algebra and
$T_{|K(A_{0})}$ is locally bounded
\endproclaim

\noindent {\it Remark.} $K(A)= K(A_{0})\oplus K(A_{1})$.

{\it Proof.} By 4.1, 4.2, and 3.6 (a), a linear transformation $t_{\pi}:
H_{\pi}\to H_{\pi}$ is defined for each irreducible $\pi$.  Suppose
$t_{\pi}$ is not bounded.  Let $x$ be an element of $K(A)$ such that
$\pi(x)\not= 0$, let $L_{x}$ be the closed left ideal generated by $x$, and
let $I_{x}$ be the closed two-sided ideal generated by $x$ (or by $L_{x}$).
Then $L_{x}\subset K(A)$ and the proof of 4.2 shows that $T(L_{x})\subset
L_{x}$ (in the notation of 4.2, $(x^{*}x)^{\frac{1}{4}}$ is in $L$).  Thus
we can apply 4.4 with $(L_{x}, I_{x})$ in the role of $(L,A)$ to obtain
that $[\pi]$ is a clopen point of $\hat I_{x}$, and hence an open point of
$\hat A$, and $\pi(I_{x})= \Cal K(H_{\pi})$.  In particular, $\pi(x)\in
\Cal K(H_{\pi})$, and since this holds for any $x$ in $K(A)$ with
$\pi(x)\not= 0$, we see that $\pi(K(A))\subset\Cal K(H_{\pi})$, and
thus $\pi(A)= \Cal K(H_{\pi})$.  This implies that $[\pi]$ is closed
in $\hat A$ and $\pi$ corresponds to an elementary direct summand of $A$.

Now let $U= \{[\pi]\in\hat A: t_{\pi}$ is not bounded\}.  By the above, $U$
is open and we now prove $U$ closed.  Assume $[\pi]$ is not in $U$ for some
irreducible $\pi$, and choose $x$ in $K(A)$ such that $\pi(x)\not= 0$.  As
above, we can cut down to $L_{x}$ and $I_{x}$, and then 4.4 implies that
$U\cap \hat I_{x}$ is finite.  Since each point of $U$ is closed, this
implies that $[\pi]$ is not in $\overline U$.

Thus $U$ induces a direct sum decomposition, $A= A_{0}\oplus A_{1}$, where
$\hat A_{1}=U$, and clearly $A_{1}$ is dual.  We replace $A$ by $A_{0}$ and
change notation.  Thus from now on we assume $t_{\pi}$ bounded for all
$\pi$.  Now we complete the proof by showing that $T_{|K(A)\cap I_{x}}$ is
bounded for all $x$ in $K(A)$.  Since $T(L_{x})\subset L_{x}$, the closed
graph theorem implies $T_{|L_{x}}$ is bounded (cf. 3.6 (d)).  Now 3.6 (b)
implies that $||t_{\pi}||\leq ||T_{|L_{x}}||$ whenever $[\pi]\in \hat
I_{x}$.  Note that $T[K(A)\cap I_{x}]\subset I_{x}$, since $Ty\in L_{y}$
for all $y$ in $K(A)$.  Then for $y$ in $K(A)\cap I_{x}, ||Ty|| =
\underset{\hat I_{x}}\to\sup ||\pi(Ty)|| = \underset{\hat I_{x}}\to\sup
||t_{\pi}\pi(y)||\leq ||T_{|L_{x}}|| \ ||y||$.  In other words,
$||T_{|K(A)\cap I_{x}}||= ||T_{|L_{x}}||$.

If $A$ is the elementary $X^{*}$-algebra $\Cal K(H),$ then $K(A) =
\Cal F(H)$, the set of bounded, finite rank operators on $H. \ \Cal
F(H)$ is spanned by (rank one) operators of the form $v\times w, \ v, w\in
H$, where $(v\times w)u = (u,w)v$.
\proclaim{Proposition 4.6} (H. Kim).  Let $A= \Cal K(H)$.

(a) The left centralizers of $K(A)$ correspond one-to-one to linear
transformations $t: H\to H$, as follows:
$$
T(v\times w) = (tv)\times w.
$$

(b) The quasi-centralizers of $K(A)$ correspond one-to-one to sesqui-linear
forms $f: H\times H\to {\Bbb C}$, as follows:
$$
\aligned
T(v_{1}\times w_{1}, v_{2}\times w_{2})&= f(v_{2},w_{1}) v_{1}\times w_{2},
\ \text{or}\\
pT(x,y)q &=f(w, x^{*}v) v\times w, \ \text{if} \ p=v\times v, \ q=
w\times w.\endaligned
$$
\endproclaim

{\it Sketch of proof.} (a) It is easy to see that for any $t$ the formula
given extends by linearity to a left centralizer.  Conversely, for any left
centralilzer $T$ and any $w$ in $H$, it is easy to see that a linear
transformation $t$ exists as in the formula.  Using $R_{u}$, for $u$  a rank
one operator, one easily shows that $t$ is independent of $w$.  (Of
course, 3.6 (a) also applies.)

(b) It is easy to see that for any $f$, the formula given extends by
bilinearity to a quasi-centralizer.  Conversely, for any quasi-centralizer
$T$ and any $v_{1},w_{2}$ in $H$, it is easy to see that a sesqui-linear
form $f$ exists as in the formula.  Using $L_{u_{1}}, R_{u_{2}}$ for rank one
operators $u_{1}, u_{2}$, one easily shows that $f$ is independent of
$v_{1},w_{2}$.

Now if $A= A_{0}\oplus A_{1}$, it is easy to see that any left centralizer
of $K(A)$ is of the form $T_{0}\oplus T_{1}$, where $T_{j}$ is a left
centralizer of $K(A_{j})$.  Also if $A_{1} =  \oplus_{i}\Cal
K(H_{i})$, a $c_{0}$-direct sum (i.e., if $A_{1}$ is dual), then $K(A_{1})=
\oplus_{i}\Cal F(H_{i})$, an algebraic direct sum (each element has
only finitely many non-zero terms).  Thus 4.5  and 4.6 (a) describe all
left centralizers of $K(A)$, but note that the decomposition, $A =
A_{0}\oplus A_{1}$, depends on the particular $T$.
\proclaim{Corollary 4.7} For any $C^{*}$-algebra $A$ the following are
equivalent:

(i) Every left centralizer of $K(A)$ is locally bounded.

(ii) $A$ has no infinite dimensional elementary direct summand.
\endproclaim
\proclaim{Theorem 4.8} If $A$ is a $C^{*}$-algebra and $T$ is a
quasi-centralizer of $K(A)$, then there is a direct sum decomposition,
$A= A_{0}\oplus A_{1}$, such that $A_{1}$ is a dual $C^{*}$-algebra and
$T_{|K(A_{0})\times K(A_{0})}$ is locally bounded.
\endproclaim

{\it Proof.} In the first three steps $\pi$ denotes a fixed irreducible
representation which does not correspond to an elementary direct summand of
$A$ and the subscript $``\pi"$ is supressed.

1) For each $x$ in $K(A), T(x,\cdot)$ is a left centralizer of $K(A)$.  By
the proof of 4.5, there is $t_{x}$ in $B(H)$ such that $\pi[T(x,y)]= t_{x}
\pi(y), \ \forall y \in K(A)$.  Since $(x,y)\mapsto T(y^{*}, x^{*})^{*}$ is
also a quasi-centralizer, the same argument gives for each $y$ in $K(A)$ a
$u_{y}$ in $B(H)$ such that $\pi[T(x,y)]= \pi(x)u_{y}, \ \forall x\in K(A)$.

2) We define a linear transformation $s: H\to H$, as follows:
$$
sv= w \ \text{if and only if} \ t_{x}v= \pi(x)w, \ \forall \ x\in K(A).
$$
Since $\pi(x)w_{1}= \pi(x)w_{2}, \ \forall x\in K(A)$, implies $w_{1}=
w_{2}, s$ is well-defined on its domain.   To see that $s$ is defined on
all of $H$, choose $y$ in $K(A)$ and $v_{0}$ in $H$ such that $\pi(y)v_{0}
= v$, for given $v$.  Then $t_{x}v= \pi[T(x,y)]v_{0}= \pi(x)u_{y}v_{0}$,
so that $sv = u_{y}v_{0}$.  It is now clear that $s$ is linear and also
that $\pi[T(x,y)]= \pi(x)s\pi(y)$.

3) We show that $s$ is bounded.  As in the proof of 4.5, choose a closed
left ideal $L$ and a closed two-sided ideal $I$ such that $L\subset K(A),
L$ generates $I$, and $\pi$ is non-trivial on $I$.  Consider the bilinear function $\pi[T(x,y)]$ for $x$ in $L^{*}$ and $y$ in $L$.  By 1) this function is
separately continuous.  It is then a well known consequence of the uniform
boundedness principle that there is a constant $c$ such that $||\pi[T(x,y)]||\leq c||x|| ||y||, \ \forall \ x\in L^{*}, \ \forall \ y\in L$.  Let $p$ be the open projection such that
$L=X(1,p)$ and note that $\pi^{**}(p)\not= 0$.  Choose a unit vector
$v_{0}$ in $\pi^{**}(p)H$.  For given vectors $v,w$ in $H$, choose, using
3.3, $x,y$ in $L$ such that $||x||= ||v||, ||y||= ||w||, \pi(x)v_{0}= v$,
and $\pi(y)v_{0}=w$.  Then since $\pi[T(x^{*},y)]= \pi(x)^{*}s\pi(y)$, we
have $|(sw,v)|= ||p_{0}\pi(x)^{*}s\pi(y)p_{0}||\leq c||v|| \ ||w||$, where
$p_{0}= v_{0}\times v_{0}$.  Thus $||s||\leq c$.

4) If $\pi$ is an irreducible which does correspond to an elementary direct
summand, then 4.6  (b) gives a sesqui-linear form $f_{\pi}: H_{\pi}\times
H_{\pi}\to {\Bbb C}.$  Let $U= \{[\pi]\in\hat A: f_{\pi}$ (as above) is not
bounded.\}.  An argument similar to part of the proof of 4.5 will show that
$U$ is closed provided we prove the following:  For each $z$ in $K(A),
U\cap\hat I_{z}$ is finite.  If this is false, there is a sequence
$\{[\pi_{n}]\}$ of distinct points in $U\cap\hat I_{z}$.  Let $I_{n}$ be
the elementary direct summand corresponding to $\pi_{n}$ and let $L_{n}=
L_{z}\cap I_{n}$.  (Each $I_{n}$ is an ideal of $A$ and $\hat I_{n}=
\{[\pi_{n}]\}$.  Although each $I_{n}$ is a direct summand of $A$, it is
not valid to write $A= A_{0}\oplus [\oplus_{n}I_{n}]$.)  For each $n$
choose $v_{n},w_{n}$ in $H_{\pi_{n}}$ such that $||v_{n}||, ||w_{n}|| <
2^{-n}$ and, $|f_{\pi_{n}}(w_{n},v_{n})| > n$.  Also choose a unit vector
$u_{n}$ in $\pi^{**}_{n}(p_{n})H_{\pi_{n}}$, where $p_{n}$ is the open
projection such that $L_{n}= X(1,p)$.  Then choose $x_{n}, y_{n}$ in
$L_{n}$ such that $||x_{n}||= ||v_{n}||, ||y_{n}||= ||w_{n}||,
\pi_{n}(x_{n})u_{n}= v_{n},$ and $\pi_{n}(y_{n})u_{n}= w_{n}$.  Then
$||T(x^{*}_{n},y_{n})|| \geq ||q_{n}T(x^{*}_{n},y_{n})q_{n}||=
|f_{\pi_{n}}(w_{n},v_{n})| > n,$
 where $q_{n}= u_{n}\times u_{n}$.  Let $x= \Sigma^{\infty}_{1}x_{n}$ and $y=
\Sigma^{\infty}_{1}y_{n}$, so that $x,y\in L_{z}\subset K(A)$.  The fact that
$T(x^{\prime},\cdot)$ is a left centralizer implies by arguments already
given that $T(x^{\prime},y^{\prime})\in L_{y^{\prime}}, \ \forall \
x^{\prime}, y^{\prime}\in  K(A)$, and similarly
$T(x^{\prime},y^{\prime})\in R_{x^{\prime}}$, the closed right ideal
generated by $x^{\prime}, \ \forall \ x^{\prime}, y^{\prime}\in K(A)$.  If
we apply this for $x^{\prime}= \Sigma_{k\not= n}x^{*}_{k}, y^{\prime} =
\Sigma_{k\not= n} y_{k}$, we see that $\pi_{n}[T(x^{*},y)]=
T(x^{*}_{n},y_{n})$, which contradicts the fact that
$||\pi_{n}[T(x^{*},y)]||\leq ||T(x^{*},y)||, \ \forall \ n$.

5)  We can now write $A= A_{0}\oplus A_{1}$, where $\hat A_{1}=U$ and
$A_{1}$ is dual.  Replace $A$ by $A_{0}$ and change notation, so that
from now on $U = \emptyset$.  Then for each irreducible $\pi$, we either have
$s_{\pi}$ in $B(H_{\pi}$) satisfying the formula in 2) or a bounded
bilinear form $f_{\pi}$ as in 4).  In the latter case there is $s_{\pi}$ in
$B(H_{\pi})$ such that $f_{\pi}(v,w) = (s_{\pi}v,w)$, so that  for all
$\pi, \pi[T(x,y)]= \pi(x)s_{\pi}\pi(y), \ \forall \ x,y\in K(A)$.

6) To complete the proof, we show that $T$ is bounded on $[K(A)\cap
I_{x}]\times [K(A)\cap I_{x}]$ for each $x$ in $K(A)$.  For any $y$ in
$K(A)$ and $\pi$ in $\hat I_{x}$, the transformation $t_{\pi}$ for the left
centralizer $T(y,\cdot)$ is $\pi(y)s_{\pi},$ which is bounded.  Therefore
$T(y,\cdot)$ is bounded on $L_{x}$ by the proof of 4.5.
Similarly, $T(\cdot,y)$ is bounded on $L^{*}_{x}$, and as above there is a
constant $c$ such that $||T(y,z)||\leq c||y|| \ ||z||, \ \forall \ y\in
L^{*}_{x}, \ \forall \ z \in L_{x}$.  As in 3) above, we then see that
$||s_{\pi}||\leq c$  whenever $[\pi]\in \hat I_{x}$.  Now arguments
already given show that $T(y,z)\in I_{x}$ whenever $x,y$ are in $K(A)\cap
I_{x}$, so that $T$ is bounded (by $c$) on $[K(A)\cap I_{x}]\times
[K(A)\cap I_{x}]$.

Again, 4.8 and 4.6 (b) describe all quasi-centralizers of $K(A)$.
\proclaim{Corollary 4.9} For any $C^{*}$-algebra $A$, the following
are equivalent:

(i) Every quasi-centralizer of $K(A)$ is locally bounded.

(ii) $A$ has no infinite dimensional elementary direct summand.
\endproclaim
\vfill\eject
\subheading{\S 5. Local linear transformations in the context of
$C^{*}$-algebras and open projections}

If $\pi: A\to B(H)$ is an irreducible representation, then there is an
isomorphism between $A^{**}$ and $B(H)\oplus M$, where $M= kernel \
\pi^{**}$.  To any rank one projection $p$ in $B(H)$ corresponds a minimal
(non-zero) projection in $A^{**}$, and all minimal projections in $A^{**}$
arise in this way.  There is a one-to-one correspondence between pure
states of $A$ and minimal projections in $A^{**}$ as follows:  Any $\varphi$
in $P(A)$ is given by $\varphi(a)= (\pi(a)v,v)$ for some irreducible $\pi$
and some unit vector $v$ in $H_{\pi}$, and $\varphi$ corresponds to the
projection $p$ with range ${\Bbb C}v$.  Then $p$ is called the {\it
support projection} of $\varphi$ because it is the smallest projection such
that $\varphi(pxp)= \varphi(x), \ \forall \ x\in A^{**}$.  As usual
we make no notational distinction between states of $A$ and the
corresponding normal states of $A^{**}$.
\proclaim{Lemma 5.1} Assume $A$ is a $C^{*}$-algebra, $p$ is a closed
projection in $A^{**}, \varphi$ is an element of $P(A)$ such that
$\varphi(p)=0$, and the support projection, $p_{0}$, of $\varphi$ is not in
$A$.  Then there is a net $\{\varphi_{i}\}$  in $P(A)$ such that
$\varphi_{i}\perp \varphi$ (i.e., $p_{i}p_{0}=0$, where $p_{i}$ is the
support projection of $\varphi_{i}), \varphi_{i}(p) = 0$, and
$\varphi_{i}\to\varphi$ weak$^{*}$.
\endproclaim

\noindent{\it Remark.} Intuitively, this lemma deals with the concept of an
isolated point, or rather a non-isolated point.

{\it Proof.} We reduce to the unital case as follows: Let $\widetilde A$ be the
result of adjoining an identity to $A$ and note that $\widetilde A^{**}\cong
A^{**}\oplus {\Bbb C}$.  Replace $p$ by $p\oplus 1$, a closed projection in
$\widetilde A^{**}$.  Since $P(A)$ can be identified with $\{\varphi\in P(\widetilde A): \varphi(0\oplus 1)=0\}$, the desired conclusion for $A$ follows
from the conclusion for $\widetilde A$.  Thus for the remainder of the
proof we assume $A$ unital.

Now consider the hereditary $C^{*}$-algebra $B$ corresponding to the
open projection $1-p_{0}$.  By Akemann's Urysohn lemma \cite{2} or \cite{3}
there is $e$ in $B$ such that $p\leq e\leq 1-p_{0}\leq 1$.  Let
$\{f_{j}\}_{j\in D_{1}}$ be an (increasing) approximate identity for $B$.
Then $\{e_{j}\}$ is also an approximate identity, where $e_{j}= e+
(1-e)^{\frac{1}{2}}f_{j}(1-e)^{\frac{1}{2}}$.  Let $D=D_{1}\times {\Bbb N}$
and for $i= (j,n)$ in $D$ choose a pure state $\varphi_{i}$ of $B$
supported by $E_{[0,\frac{1}{n}]}(e_{j})$, where the last symbol denotes a
spectral projection (in $B^{**}$) of $e_{j}$.  If this were impossible,
i.e., if the closed projection $E_{[0,\frac{1}{n}]}(e_{j})$ were zero, then
$e_{j}$ would be invertible in $B$ and $B$ would be a unital
$C^{*}$-algebra.  This would imply $p_{0}\in A$, a contradiction.  Since
$p\leq E_{\{1\}}(e_{j}), \varphi_{i}(p) = 0$ for $n > 1$, and clearly
$\varphi_{i}\perp\varphi$.  We see as follows that $\varphi_{i}\to 0$ in
the weak$^{*}$ topology of $B^{*}$: Given $b\in B, \ b\geq 0, \ b\approx
e_{j}^{\frac{1}{2}}be_{j}^{\frac{1}{2}}$ for $j$ sufficiently large.  Hence
$\varphi_{j,n}(b)\approx
\varphi_{j,n}(e^{\frac{1}{2}}_{j}be_{j}^{\frac{1}{2}})\leq
||b||\varphi_{j,n}(e_{j})\leq \frac{1}{n}||b||$.  This means that any
weak$^{*}$ cluster point of $\{\varphi_{i}\}$ in $A^{*}$ is supported by
$p_{0}$ and hence $\varphi_{i}\to\varphi$.
\proclaim{Theorem 5.2} Assume $A$ is a $C^{*}$-algebra, $p_{1}$ and $p_{2}$
are open projections in $A^{**}$, and $T: X(p_{1},p_{2})\to X(1,p_{2})$ is
a local linear transformation.  Then $T$ is purely decomposable.
\endproclaim

{\it Proof.} It is sufficient to show $\varphi(x^{*}x)=0$ implies
$\varphi[(Tx)^{*}(Tx)]=0$ for $\varphi$ a pure state supported by $p_{2}$
and $x$ in $X(p_{1},p_{2})$.  Let $p_{0}$ be the support projection of
$\varphi$, so that $p_{0}\leq p_{2}$.  Since for all $x, \varphi(x^{*}x)=0
\Leftrightarrow xp_{0}=0$, the result to be proved is immediate if $p_{0}$
is open.  Therefore we may assume $p_{0}\notin A$.  

Suppose then that $\varphi(x^{*}x)=0$ and $\varphi[(Tx)^{*}(Tx)]> 1$.  We
will construct recursively pure states $\varphi_{n}$, open projections
$q_{n}$, and elements $y_{n}$ of $X(p_{1},p_{2})$ such that $\overline
q_{n}\perp\overline q_{m}$ for $n\not= m, \ \overline q_{n}\perp p_{0}, \
\overline q_{n}\leq p_{2}, \ \varphi_{n}$  is supported by $q_{n}, \
y_{n}q_{k}=0$ for $k < n, \ ||y_{n}|| < 2^{-n}$, and $||T(y_{1}+\dots +
y_{n})w_{\varphi_{n}}|| > n$.  Here, for any projection $q$ in $A^{**},
\overline q$ denotes the smallest closed projection majorizing $q$ (cf.
\cite{1}), and for any state $\psi$ of $A, \ w_{\varphi}$ denotes the cyclic
vector in the GNS construction, so that $||yw_{\psi}|| =
[\psi(y^{*}y)]^{\frac{1}{2}}$ for any $y$ in $A^{**}$.  Then if $y=
\Sigma^{\infty}_{1}y_{n}, \ yq_{n}= (y_{1}+\dots + y_{n})q_{n}$.  Hence by
locality $||(Ty)w_{\varphi_{n}}||= ||T(y_{1}+\dots +
y_{n})w_{\varphi_{n}}|| > n$.  This is a contradiction, since
$||(Ty)w_{\varphi_{n}}|| \leq ||Ty||, \ \forall \ n$.

Assume we already have $\varphi_{1},\dots , \varphi_{n-1}, q_{1},\dots,
q_{n-1}$, and $y_{1},\dots, y_{n-1} \ (n\geq 1)$.  Choose $\lambda >
||T(y_{1}+\dots + y_{n-1})|| + n \ (\lambda > 1$ if $n=1$).  Let $\delta=
\lambda^{-1}2^{-n}$ and choose a continuous function $f_{\delta}:
[0,\infty)\to [0,1]$ such that $f_{\delta}=1$ in a neighborhood of $0$ and
$f_{\delta}(t)=0$ for $t\geq\delta$.  Let $x^{\prime}= xf_{\delta}(|x|)$,
where $|x|=(x^{*}x)^{\frac{1}{2}}$.  Then $||x^{\prime}|| < \delta$ and
$(x^{\prime}-x)q=0$ where $q$ is an open projection of the  form
$E_{[0,\epsilon)}(|x|)$.  Since $p_{0}\leq q$, it follows from locality that
$[T(x^{\prime}-x)]p_{0}=0$.  Thus, $||(Tx^{\prime}) w_{\varphi}||=
||(Tx)w_{\varphi}|| > 1$.  Now by Akemann's Urysohn lemma choose $h$ in $A$
such that $p_{0}\leq h\leq p_{2}- (\overline q_{1}+\dots + \overline
q_{n-1})$.  The fact that the projection on the  right is open, i.e.,
$1-p_{2} + \overline q_{1}+\dots + \overline q_{n-1}$ is closed, follows
from \cite{1}.   Let $g: [0,1]\to [0,1]$ be a continuous function such that
$g(0)=0$ and $g=1$ in a neighborhood of $1$.  Let $y_{n}= \lambda
x^{\prime}g(h)$.  Then $||y_{n}|| < 2^{-n}$, and by locality
$||(Ty_{n})w_{\varphi}||= \lambda||(Tx^{\prime})w_{\varphi}|| > \lambda$.
Now use 5.1 to find $\varphi_{n}$ in $P(A)$ such that
$\varphi_{n}\perp\varphi, \ \varphi_{n}(1-p_{2}+ \overline q_{1}+\dots +
\overline q_{n-1})=0$ and $\varphi_{n}[(Ty_{n})^{*}(Ty_{n})] >\lambda^{2}$.
Therefore $||(Ty_{n})w_{\varphi_{n}}|| > \lambda$ and hence $||T(y_{1}+\dots
+ y_{n})w_{\varphi_{n}}|| > n$.  Finally, if $r_{n}$ is the support
projection of $\varphi_{n}$, then we can find, by Akemann's Urysohn
lemma, $h^{\prime}$ in $A$ such that $r_{n}\leq h^{\prime}\leq p_{2}-
(\overline q_{1}+\dots + \overline q_{n-1})-p_{0}$.  Then let
$q_{n}=E_{(\frac{1}{2},1]}(h^{\prime})$, so that $\overline q_{n}\leq
E_{[\frac{1}{2},1]}(h^{\prime})$.  This completes the recursive
construction and the proof of the theorem.
\proclaim{Corollary 5.3} Under the same hypotheses the ideal of $A$
generated by $X(p_{1},p_{2})$ can be written as a direct sum, $A_{0}\oplus
A_{1}\oplus\dots \oplus A_{n}$, such that $X(p_{1},p_{2})=
\oplus^{n}_{0}[X(p_{1},p_{2})\cap A_{i}], A_{i}= \Cal K(H_{i})$ for $i >
0$, and $T(x_{0}\oplus x_{1}\oplus \dots\oplus x_{n})= t_{0}x_{0}\oplus
t_{1}x_{1}\oplus\dots\oplus t_{n}x_{n}$ for some $t_{0}$ in $A^{**}_{0}$
and (possibly discontinuous) linear transformations $t_{i}:
\pi^{**}_{i}(p_{1})H_{i}\to H_{i}$, where $\pi_{i}$ is the irreducible
representation corresponding to $A_{i}$.
\endproclaim

{\it Proof.} There is no loss of generality in assuming $X(p_{1},p_{2})$
generates $A$ as an ideal.  Note that if $p_{1}=1$, then $X(p_{1},p_{2})$
is a left ideal, and we are in the situation of 4.3 and 4.4.  The main
point is to generalize 4.3 and 4.4, which will give us the desired direct
sum decomposition so that $T_{|X(p_{1},p_{2})\cap A_{0}}$ is bounded, and
$t_{i}= t_{\pi_{i}}$ for $i > 0$.

As already remarked, the only parts of the proofs of 4.3 and 4.4 that
require change are the justifications of cutting down to an ideal $I$ of
$A$.  The fact that $X(p_{1},p_{2})\cap I$ generates $I$ as an ideal can be
proved as follows (it requires no proof for readers familiar with the
theory of imprimitivity bimodules): Let $B_{i}= her(p_{i})$, so that
$X(p_{1},p_{2})= [B_{1}AB_{2}]^{-}$, and let $Y= X(p_{1},p_{2})\cap I$.
Then $[IYI]^{-}\supset [I(B_{1}IB_{2})I]^{-}\supset
[(IA)B_{1}(AIA)B_{2}(AI)]^{-}\supset I(AB_{1}A)^{-}I (AB_{2}A)^{-}I=
IAIAI=I$.  Here we used the fact that $B_{i}$ generates $A$ as an ideal,
since $X(p_{1},p_{2})$ is clearly contained in the ideal generated by
$B_{i}$.  The fact that $T(Y)\subset X(1,p_{2})\cap I$ follows from pure
decomposability.  We must show that for an irreducible $\pi$ which is
non-trivial on $I, t_{\pi}$ is the same whether computed for $A$ and
$X(p_{1},p_{2})$ or $I$ and $Y$.  Since the first version of $t_{\pi}$ is
an extension of the second, we need only show the domains are the same.
Now $Y= X(zp_{1},zp_{2})$, where $z$ is the open central projection
corresponding to $I$; i.e., $I = \{a\in A: za=a\}$.If $\pi$ is non-trivial
on $I$, then $\pi^{**}(z)= 1$ and hence $\pi^{**}(zp_{1})= \pi^{**}(p_{1})$.  

Now we have $A= A_{0}\oplus\dots\oplus A_{n}$ and the fact that
$X(p_{1},p_{2})= \oplus^{n}_{0}[X(p_{1},p_{2})\cap A_{i}]$ is immediate
from the fact that $p_{i}$ can be written, $p_{i}= \oplus^{n}_{j=0}p_{ij}$
with $p_{ij}$ in $A^{**}_{j}$.  Since $T_{|X(p_{1},p_{2})\cap A_{0}}$ is
bounded, 3.10 yields $t_{0}$ in $A^{**}_{0}$ with the desired property.

\noindent{\it Remark.} Of course, $t_{0}$ is not an arbitrary element of
$A^{**}_{0}$, since necessarily $t_{0}[X(p_{1},p_{2})\cap A_{0}]\subset
A_{0}$.  One should regard $t_{0}$ as a kind of left multiplier (this
interpretation possibly should be used only in the special case where
$T(X(p_{1},p_{2}))\subset X(p_{1},p_{2}))$.  It has already been mentioned
that $X(p_{1},p_{2})$ is a Hilbert $B_{1}-B_{2}$  bimodule (an
imprimitivity bimodule, even, if $X(p_{1},p_{2})$ generates $A$ as an
ideal), and Hilbert bimodules have many properties in common with
$C^{*}$-algebras.  However, every left centralizer of a $C^{*}$-algebra is
bounded by \cite{22, 3.12.2}, and in certain cases left centralizers of
$X(p_{1},p_{2})$, of the sort considered in 5.3, are not bounded.  The
exceptions occur only when the ideal generated by $X(p_{1},p_{2})$ has an
elementary direct summand $\Cal K(H_{\pi})$ such that $\pi^{**}(p_{2})$
has finite rank and $\pi^{**}(p_{1})$ has infinite rank.
\proclaim{Corollary 5.4} If $p_{1}=p_{2}$, then $T$ is bounded.  In
particular, any local linear transformation from $A$ to itself is bounded.
\endproclaim
\noindent {\it Remark.} Since $T$ is given only as a local linear
transformation, not as a left centralizer, the last sentence of 5.4 is a
new result.  Of course, {\it a postiori}, $T$ is a left centralizer.
\proclaim{Proposition 5.5} Assume $A$ is a $C^{*}$-algebra, $p_{1}$ and
$p_{2}$ are open projections in $A^{**}$, and $T$ is a linear
transformation from $X(p_{1},p_{2})$ to $X(1,p_{2})$.  Then the following
are equivalent:

(i) $T(xb)= (Tx)b,$ for all $b$ in $her(p_{2})$, the hereditary
$C^{*}$-subalgebra of $A$ supported by $p_{2}$, and all $x$ in
$X(p_{1},p_{2})$.

(ii) $T$ is decomposable.

(iii) $T$ is local.

(iv) $T$ is purely decomposable.
\endproclaim

{\it Proof.} (i) $\Rightarrow$ (ii).  We can write $x= y|x|^{\frac{1}{2}}$
for some $y$ in $X(p_{1},p_{2})$ (cf. \cite{22, 1.4.5}).  Then
$|x|^{\frac{1}{2}}\in her(p_{2})$, and hence $Tx= (Ty)|x|^{\frac{1}{2}}$.
Since $\In \ x = \In \ |x|^{\frac{1}{2}}$, this shows that $\In \ Tx\leq \In
\ x$.

(ii) $\Rightarrow$ (iii) is trivial.

(iii) $\Rightarrow$ (iv) is 5.2.

(iv) $\Rightarrow$ (i).  If $\pi$ is irreducible and $x,b$ are as in (i),
then by 3.6 (a), $\pi[T(xb)]= t_{\pi}\pi(xb)$ and $\pi[(Tx)b] =
\pi(Tx)\pi(b)= t_{\pi}\pi(x)\pi(b) = t_{\pi}\pi(xb)$ (in case
$\pi^{**}(p_{1})$ or $\pi^{**}(p_{2})$ is $0$, both above are $0$).  Since
there are enough irreducible representations to distinguish elements of $A$
(i.e., $A$ is of type $\Cal U$), this implies $T(xb)= (Tx)b$.

\noindent {\it Remark 5.6.} It is possible to combine the subjects of
Sections 4 and 5 by considering local operators $T$ defined on
$X(p_{1},p_{2})\cap K(A)$.  The analogues of 5.2 and 5.5 are true in this
context and also the analogue of 5.3 and 4.5: The ideal of $A$ generated by
$X(p_{1},p_{2})$ can be written as $A_{0}\oplus A_{1}$ such that $A_{1}$ is
dual and $T_{|X(p_{1},p_{2})\cap K(A)\cap A_{0}}$ is locally bounded.  The
proofs do not require much in the way of new arguments:

1. In the proof of 5.2, $y$ is in the closed right ideal generated by $x$
and hence $x$ in $K(A)$ implies $y$ in $K(A)$.

2. For $x$ in $K(A)\cap X(p_{1},p_{2})$, let $q$ be the support projection
of $x^{*}x$.  Then $q$ is open, $q\leq p_{2}$, and $x\in X(p_{1},q)\subset
K(A)$.  This provides a substitute for $L_{x}$, as used in the proof of
4.5.  $(L_{x}= X(1,q).)$

3. If $I$ is an ideal of $A$, then $K(A)\cap I$ may be strictly larger then
$K(I)$.  Thus if we have a direct sum decomposition, $A_{0}\oplus A_{1}$,
as above, it is not obvious how to describe the domain of $T$ in terms of
$A_{0}$ and $A_{1}$.  Nevertheless, the fact that $K(A)$ is a union of
closed left ideals can be used to show that the domain of $T$ is compatible
with the direct sum decomposition.  Also, writing $A_{1}= \underset
i\to\oplus \Cal K(H_{i})$ and using the same idea, we can show that
every element of $K(A)\cap X(p_{1},p_{2})\cap A_{1}$ has only finitely many
non-zero components (provided each $t_{i}$ is discontinuous).  If this were
false a construction like the one in the last paragraph of the proof of 4.3
could be carried out (within the set $X(p_{1},q)$ of point 2 above) to give
a contradiction.  (Thus, {\it a postiori}, $K(A)\cap X(p_{1},p_{2})\cap
A_{1}\subset K(A_{1}).$
\vfill\eject
\bye